\documentclass[reqno]{amsart}

\usepackage{amssymb,amsmath,amsthm,amsfonts,fancyhdr}
\usepackage{setspace}
\usepackage{cite}

\newtheorem{thm}{Theorem}[section]
\newtheorem{defn}{Definition}[section]
\newtheorem{lemma}{Lemma}[section]

\newtheorem{re}{Remark}[section]

\numberwithin{equation}{section}

\begin{document}

\title[Remarks on area maximizing hypersurfaces]{Remarks on area maximizing hypersurfaces over $\mathbb{R}^n\backslash\{0\}$ and exterior domains}

\author{Guanghao Hong}

\address{School of Mathematics and Statistics, Xi'an Jiaotong University, Xi'an, P.R.China 710049.}

\email{ghhongmath@xjtu.edu.cn}

\begin{abstract}
In this note, we provide a complete classification for entire area maximizing hypersurfaces having an isolated singularity. We also construct an interesting illustrated example. For area maximizing hypersurfaces over exterior domains, we obtain a partial result on their asymptotic behavior at infinity. We also establish the solvability of exterior Dirichlet problems for area maximizing hypersurfaces.
\end{abstract}


\maketitle

\section{Introduction}

Area maximizing hypersurfaces are the graphs of the solutions $u$ to the variational problem for the area functional $$\int_{\Omega} \sqrt{1-|Du|^2}$$ in Lorentz-Minkowski space $\mathbb{L}^{n+1}$. Maximal hypersurfaces are the graphs of the smooth solutions $u$ to the corresponding Euler-Lagrange equation (known as maximal surface equation)
\begin{equation}
 div(\frac{D u}{\sqrt{1-|D u|^2}})=0, \ \ |Du|<1.
\end{equation}
Area maximizing hypersurfaces are weaker than maximal hypersurfaces.

Calabi [Ca68] ($n\leq 4$) and Cheng-Yau [CY76] (all dimensions) showed that the only entire maximal hypersurfaces are spacelike hyperplanes. Bartnik-Simon [BS82] proved existence of solutions for the variational problem and for the equation with appropriate prescribed boundary values on bounded domains and compared these two kinds of solutions. Base on the results in [BS82], Bartnik observed that the only entire area maximizing hypersurfaces are weakly spacelike hyperplanes.

Ecker [Ec86] studied isolated singularity problem for area maximizing hypersurfaces. He achieved two main results. The first one is that ``Isolated singularities of area maximizing hypersurfaces in $\mathbb{L}^{n+1}$ are light-cone-like.'' The second one is that ``Entire maximal hypersurfaces in $\mathbb{L}^{n+1}$ having an isolated singularity are the Lorentz transformations of the known radially symmetric (about the singularity) solutions.'' However, the only thing he said about entire area maximizing hypersurfaces with an isolated singularity is that ``If such a hypersurface is a cone, then it is the upper or lower light cone.'' In this note, we provide a complete classification for entire area maximizing hypersurfaces with an isolated singularity. See the precise notations and definitions in Section 2.

\begin{thm}
Let $M$ (gragh of a weakly spacelike function $u$) be an entire area maximizing hypersurface with an isolated singularity at $0\in \mathbb{L}^{n+1}$, then exact one of the four situations happens.

(i) $M$ is an entire maximal hypersurface with an isolated singularity;

(ii) $M$ is the upper light cone, \textit{i.e.}, $u(x)=|x|$;

(iii) $M$ is the lower light cone, \textit{i.e.}, $u(x)=-|x|$;

(iv) $M$ is asymptotic to a hyperplane with slope 1 at infinity, \textit{i.e.}, $$u(x)=a\cdot x+o(|x|)\ \ \mbox{as}\ |x|\rightarrow \infty$$ for some unit vector $a\in \mathbb{R}^n$. Moreover, in this case, either $u(x)\leq a\cdot x$ and $u(ta)=t$ for $t\leq 0$ or $u(x)\geq a\cdot x$ $u(ta)=t$ for $t\geq 0$.
\end{thm}

In Section 3.2, we construct an area maximizing hypersurface satisfying (iv) of Theorem 1.1. The construction is copied from [SWY08] (Page 4). Their construction is for infinity harmonic functions.

In [HY18], we achieved a complete description of asymptotic behavior at infinity for exterior maximal hypersurfaces. However, for exterior area maximizing hypersurfaces, we only have the following partial result so far.

\begin{thm}

Let $M$ (gragh of a weakly spacelike function $u\in C(\mathbb{R}^n\backslash A^o)$, where $A$ is a bounded closed set and $A^o$ is its interior) be an area maximizing hypersurface on $\mathbb{R}^n\backslash A$. We assume $u$ (or $M$) is not a maximal hypersurface on the whole $\mathbb{R}^n\backslash \mbox{Conv}(A)$, where $\mbox{Conv}(A)$ denotes the convex hull of $A$. Then exact one of the three situations happens.

(i) $u(x)-|x-x_0|$ attains its maximum and minimum on $\partial A$ for any $x_0\in A$;

(ii) $u(x)+|x-x_0|$ attains its maximum and minimum on $\partial A$ for any $x_0\in A$;

(iii) There exists a unit vector $a\in \mathbb{R}^n$ such that all blowdowns of $u(x)$ have the form $a\cdot x+o(|x|)$ (as $|x|\rightarrow \infty$) and one of these blowdowns is $a\cdot x$. Moreover, $u(x)-a\cdot x$ is bounded from one side.
\end{thm}

One certainly wish to improve the formulation in (iii) of Theorem 1.2 to $$u(x)=a\cdot x+o(|x|)\ \ \mbox{as}\ |x|\rightarrow \infty.$$ But this is not easy since the solution $u$ is real degenerate in this situation.

Finally, we establish the solvability of exterior Dirichlet problem for area maximizing hypersurfaces. Given a bounded closed set $A\subset \mathbb{R}^n$, we say a boundary value function $g:\partial A \rightarrow \mathbb{R}$ is admissible if $g$ is bounded and there exists a weakly spacelike function $\psi$ in $\mathbb{R}^n\backslash A$ such that $\psi=g$ on $\partial A$ (in the sense of (1.1) in [BS82]).

\begin{thm}
Let $A\subset \mathbb{R}^n$ be a bounded closed set and $g:\partial A \rightarrow \mathbb{R}$ be admissible. Then we have the followings.

(i) For any $x_0\in A$, there exists an area maximizing function $u$ on $\mathbb{R}^n\backslash A$ satisfying $u=g$ on $\partial A$ and $u(x)-|x-x_0|$ attains its maximum and minimum on $\partial A$;

(ii) For any $x_0\in A$, there exists an area maximizing function $u$ on $\mathbb{R}^n\backslash A$ satisfying $u=g$ on $\partial A$ and $u(x)+|x-x_0|$ attains its maximum and minimum on $\partial A$;

(iii) Given any unit vector $a\in \mathbb{R}^n$, there exists an area maximizing function $u$ on $\mathbb{R}^n\backslash A$ satisfying $u=g$ on $\partial A$ and $u(x)-a\cdot x$ attains its maximum and minimum on $\partial A$.
\end{thm}

For maximal hypersurfaces, we completely established the existence and uniqueness of solutions to exterior Dirichlet problems [HY18]. But for area maximizing hypersurfaces, we know almost nothing about the uniqueness of the solution.

Area maximizing hypersurfaces and infinity harmonic functions appear to enjoy significant similarities. We encourage the readers to compare [SWY08] and our another recent paper [HZ18] with [Ec86] and this work.

The paper is organized as follows. In Section 2, we set up some notations and state some results from [BS82] and [Ec86] for our later use. In Sections 3, we prove Theorem 1.1 and construct the above-mentioned example. In section 4, we prove Theorems 1.2 and 1.3.

\section{Preliminaries}

 We denote the Lorentz-Minkowski space by $\mathbb{L}^{n+1}=\{X=(x,t): x\in \mathbb{R}^n, t\in \mathbb{R}\}$, with the metric $\sum_{i=1}^n dx_i^2-dt^2$. And $\langle \cdot,\cdot\rangle$ denotes the inner product in $\mathbb{L}^{n+1}$ with the signs $(+,\cdots,+,-)$. The light cone at $X_0=(x_0,t_0)\in\mathbb{L}^{n+1}$ is defined by $$C_{X_0}=\{X\in\mathbb{L}^{n+1}:\langle X-X_0,X-X_0\rangle=0\}.$$ The upper and lower light cones will be denoted by $C^+_{X_0}$ and $C^-_{X_0}$ respectively.

 Let $M$ be an $n$-dimensional hypersurface in $\mathbb{L}^{n+1}$ which can be represented as the graph of $u\in C^{0,1}(\Omega)$, where $\Omega$ is a open set in $\mathbb{R}^n$. We say $M$ (or $u$) is \textit{weakly spacelike} if $|Du|\leq 1$ a.e. in $\Omega$.

\begin{thm}[see \cite{BS82}]
Let $\Omega\subset \mathbb{R}^n$ be a bounded domain and let $\varphi: \partial \Omega\rightarrow R$ be a bounded function. Then the variational problem
\begin{equation}
\int_{\Omega}\sqrt{1-|Dv|^2}\rightarrow \mbox{max}\ \     \mbox{in}\   K
\end{equation}
where $K=\{v\in C^{0,1}(\Omega): |Dv|\leq 1$ a.e. in $\Omega$, $v=\varphi$ on $\partial\Omega\}$
has a unique solution $u$ if and only if the set $K$ is nonempty.
\end{thm}

\begin{defn}
A weakly spacelike function $u\in C(\Omega)$ ($\Omega\subset \mathbb{R}^n$ is not necessarily bounded) is called \textit{area maximizing} if it solves the variational problem (2.1) with respect to its own boundary values for every bounded subdomain in $\Omega$. The graph of $u$ is called an \textit{area maximizing hypersurface}.
\end{defn}

\begin{lemma}[see \cite{BS82}]
If $\{u_k\}$ is a sequence of area maximizing functions in $\Omega$ and $u_k\rightarrow u$ in $\Omega$ locally uniformly, then $u$ is also an area maximizing function.
\end{lemma}

 One key result in [BS82] (Theorem 3.2) is that if an area maximizing hypersurface contains a segment of light ray, then it contains the whole of the ray extended all the way to the boundary or to infinity.

\begin{thm}[see \cite{BS82}]
The solution $u$ of (2.1) is smooth and solves equation (1.1) in $$reg\ u:=\Omega\backslash sing\ u$$ where $$sing\ u:=\{\overline{xy}:x,y\in\partial\Omega, x\neq y,\overline{xy}\subset\Omega \ \mbox{and}\ |\varphi(x)-\varphi(y)|=|x-y|\}.$$ Furthermore $$u(tx+(1-t)y)=t\varphi(x)+(1-t)\varphi(y),\ 0<t<1$$ where $x,y\in\partial\Omega$ are such that $\overline{xy}\subset\Omega$ and $|\varphi(x)-\varphi(y)|=|x-y|$.
\end{thm}

\begin{thm}[Bartnik, see Theorem F in \cite{Ec86}]
Entire area maximizing hypersurfaces in $L^{n+1}$ are weakly spacelike hyperplanes.
\end{thm}

\begin{defn}
A weakly spacelike hypersurface $M$ in $L^{n+1}$ containing $0$ is called an area maximizing hypersurface having an isolated singularity at $0$ if $M\backslash \{0\}$ is area maximizing but $M$ cannot be extended as an area maximizing hypersurface into $0$.
\end{defn}

Ecker proved that the isolated singularities of area maximizing hypersurface are light cone like (Theorem 1.5 in [Ec86]).

For a weakly spacelike exterior hypersurface $M$ (\textit{i.e.}, $u$ is defined on an exterior domain $\mathbb{R}^n\backslash A$ with $A$ bounded), we define $M_r=r^{-1}M$ with $r>0$ is the graph of $u_r(x)=r^{-1}(rx)$. If for some $r_j\rightarrow +\infty$, $u_{r_j}(x)$ converge locally uniformly to a function $u_{\infty}(x)$ on $\mathbb{R}^n\backslash \{0\}$, then $u_{\infty}$ (its graph $M_{\infty}$) is called a blowdown of $u$ ($M$). Note that by weakly spacelikeness, Arzela-Ascoli theorem always ensures the existence of blowdowns. By Lemma 2.1, $u_{\infty}(x)$ ($M_{\infty}$) is area maximizing on $\mathbb{R}^n\backslash \{0\}$ and $u_{\infty}(0)=0$. The following lemma is especially useful.

\begin{lemma}[Lemma 1.10 in \cite{Ec86}]
Let $M$ be an entire area maximizing hypersurface having an isolated sigularity at $0$ and assume that some blowdown of $M$ also has an isolated singularity at $0$. Then $M$ has to be either $C^+_0$ or $C^-_0$.
\end{lemma}

\section{Entire area maximizing hypersurfaces with an isolated singularity}

\subsection{Asymptotic behavior}

\begin{proof}[Proof of Theorem 1.1]
Suppose $M$ is not an maximal hypersurface on the whole $\mathbb{R}^n\backslash\{0\}$. By Theorem 2.2, $M$ contains a (past-directed or further-directed) light ray emanating from $0\in \mathbb{L}^{n+1}$. We assume this light ray is past-directed, that is, for some unit vector $a\in \mathbb{R}^n$, $u(ta)=t$ for $t\in (-\infty,0]$. There are two alternative situations.

(a) Some blowdown of $M$ has an isolated singularity at $0$. In this case, By Lemma 2.2,  $M$ is $C^+_0$ or $C^-_0$. Since we assumed that $M$ contains a past-directed light ray, $M=C^-_0$.

(b) All blowdowns of $M$ are entire area maximizing hypersurfaces. By Theorem 2.3, These blowdowns are weakly spacelike hyperplanes. Obviously, all blowdowns contain the light ray $\{(ta,t): t\leq 0\}$. The only choice for these hyperplanes is the graph of $u_{\infty}(x)=a\cdot x$. That is to say, the blowdown of $u$ is unique and it is $u_{\infty}(x)=a\cdot x$. Thus we have $$u(x)=a\cdot x+o(|x|)\ \ \mbox{as}\ |x|\rightarrow \infty.$$ We claim one more thing: $$u(x)\leq a\cdot x\ \  \mbox{in}\ \mathbb{R}^n.$$ For simplicity, we assume $a=e_n$. Since $M$ is weakly spacelike,
\begin{eqnarray*}
u(x) & \leq & u(te_n)+|x-te_n|\ \ \  (\mbox{for}\  t\leq 0) \\
& = & t+ \sqrt{|x'|^2+(x_n-t)^2} \rightarrow x_n\ \ \ (\mbox{as}\  t\rightarrow-\infty).
\end{eqnarray*}

We have proved that in case that $M$ contains a past-directed light ray (say $\{(ta,t): t\leq 0\}$), we have that either $u(x)=-|x|$ or $u(x)=a\cdot x+o(|x|)$ and $u(x)\leq a\cdot x$. Similarly, in case that $M$ contains a future-directed light ray (say $\{(ta,t): t\geq 0\}$), we have that either $u(x)=|x|$ or $u(x)=a\cdot x+o(|x|)$ and $u(x)\geq a\cdot x$.

\end{proof}

\subsection{An example}
\

\vspace{0.618ex}
For any $k=2,3,\cdots$ and any $\theta \in [0,1]$, let $w^{\theta}_k(x)$ be the solution to (2.1) in $B_k\backslash \{0\}$ with the boundary condition $w(0)=0$ and $w(x)=\theta x_n-(1-\theta)k$ on $\partial B_k$. For each $k$, it is easy to see that $w^{0}_k(x)=-|x|$ in $B_k$ (implying $w^{0}_k(e_n)=-1$) and $w^{1}_k(x)=x_n$ (implying $w^{1}_k(e_n)=1$). By continuity, there exists some $\theta(k)\in (0,1)$ satisfying $w^{\theta(k)}_k(e_n)=0$. Consider the sequence of functions $\{w^{\theta(k)}_k(x)\}$. Due to the weakly spacelikeness and $\{w^{\theta}_k(0)=0\}$, the family of functions $\{w^{\theta(k)}_k(x)\}$ is uniformly bounded and equicontinuous on any compact set $K\subset\subset \mathbb{R}^n$. By Arzela-Ascoli theorem, up to a subsequence, $$w^{\theta(k)}_k(x)\rightarrow w(x) \ \ \mbox{locally uniformly in}\ \mathbb{R}^n.$$ By Lemma 2.1, $w(x)$ is area maximizing in $\mathbb{R}^n\backslash \{0\}$. From the construction, we can see that $w(te_n)=t$ for $t\in (-\infty,0]$ (because $w^{\theta}_k(te_n)=t$ for $t\in (-k,0]$) and $w(e_n)=0$. By Theorem 1.1, we know $$w(x)=x_n+o(|x|)\ \ \mbox{as}\ |x|\rightarrow \infty$$  and $w(x)\leq  x_n$ in $\mathbb{R}^n$.

We claim one more thing: $w$ is smooth and solves (1.1) in $\mathbb{R}^n\backslash \{te_n:t\in (-\infty,0]\}$. If this is not true, by Theorem 2.2, the graph of $w$ contains another light ray (other then $\{(te_n,t): t\leq 0\}$). Because $w(x)=x_n+o(|x|)$, this new light ray can only be $\{(te_n,t): t\geq 0\}$. But this is also impossible since $w(e_n)=0\neq 1$.

\section{Area maximizing hypersurfaces over exterior domains}

\subsection{Asymptotic behavior}

\begin{proof}[Proof of Theorem 1.2]
Since we assumed $M$ is not a maximal hypersurface on the whole $\mathbb{R}^n\backslash \mbox{Conv}(A)$, by Theorem 2.2, $M$ contains a light ray emanating from a boundary point of $A$. We assume this light ray is past-directed. For simplicity, we assume this boundary point is $0$ and this light ray is $\{(te_n,t): t\leq 0\}$. We also assume $A\subset \bar{B}_1$ without loss of generality. There are two alternative situations.

(a) There exists a blowdown $u_{\infty}$ of $u$, such that some blowdown of $u_{\infty}$ has an isolated singularity at $0$. In this case, By Lemma 2.2,  $u_{\infty}(x)=|x|$ or $-|x|$. Since $u_{\infty}(te_n)=t$ for $t\leq 0$ (because $u(te_n)=t$ for $t\leq 0$), we have $u_{\infty}(x)=-|x|$.

For any $x_0\in A$, denote $\max\limits_{\partial A}(u(x)+|x-x_0|):=c^+$ and $\min\limits_{\partial A}(u(x)+|x-x_0|):=c^-$. Then $$c^- -|x-x_0|\leq u(x)\leq c^+ -|x-x_0| \ \ \mbox{on}\  \partial A.$$ By weakly spacelikeness of $u$, we immediately have $$c^- -|x-x_0|\leq u(x)\ \ \mbox{in}\ \mathbb{R}^n\backslash A.$$ We need to show $$u(x)\leq c^+ -|x-x_0|\ \ \mbox{in}\ \mathbb{R}^n\backslash A.$$

Suppose $$u_{\infty}(x)=\lim_{k\rightarrow \infty}\frac{u(r_kx)}{r_k}\ \ \mbox{in}\ \mathbb{R}^n\backslash \{0\}$$ for a sequence of $r_k\rightarrow \infty$.
For any $\epsilon>0$, there is $\bar{k}$ such that $$u(x)\leq c^+-(1-\epsilon)|x-x_0| \ \  \mbox{on} \ \partial B_{r_k}$$ for all $k\geq \bar{k}$. Note that $c^+-(1-\epsilon)|x-x_0|$ is a supersolution to (1.1) in $\mathbb{R}^n\backslash \{x_0\}$, by comparison principle (see [BS82]), we have $$u(x)\leq c^+ -(1-\epsilon)|x-x_0|\ \ \mbox{in}\ (B_{r_{\bar{k}}}\backslash A)\cup (\cup_{j=1}^{\infty}B_{r_{\bar{k}+j}}\backslash B_{r_{\bar{k}+j-1}})=\mathbb{R}^n\backslash A.$$ Letting $\epsilon\rightarrow 0$, we have $$u(x)\leq c^+ -|x-x_0|\ \ \mbox{in}\ \mathbb{R}^n\backslash A.$$

(b) For any blowdown $u_{\infty}$ of $u$, all blowdowns of $u_{\infty}$ are entire area maximizing hypersurfaces. By Theorem 2.3, These blowdowns are weakly spacelike hyperplanes. Since all blowdowns contain the light ray $\{(te_n,t): t\leq 0\}$, all these hyperplanes have to be $\{(x,x_n): x\in \mathbb{R}^n\}$. That is to say, the blowdown of $u_{\infty}$ is unique and it is the function $V(x)=x_n$. Thus we have $$u_{\infty}(x)=x_n+o(|x|)\ \ \mbox{as}\ |x|\rightarrow \infty.$$ We also have $$u_{\infty}(x)\leq x_n\ \  \mbox{in}\ \mathbb{R}^n.$$

Next we prove that one of the blowdowns of $u$ is the function $V(x)=x_n$. We denote $m:=\max_{\partial A}u$, then $u\leq m+|x|$ by weakly spacelikeness of $u$.
We claim that there exists a sequence of points $x_k\in \mathbb{R}^n\backslash B_2$ with $r_k:=|x_k|\rightarrow \infty$ satisfying $u(x_k)>m+(1-\frac{1}{k})|x_k|$. Because otherwise, $$u\leq m+(1-\epsilon)|x|\ \ \mbox{in}\ \mathbb{R}^n\backslash B_R$$ for some $R>2$ and $\epsilon>0$. That implies any blowdown $u_{\infty}\leq (1-\epsilon)|x|$, contradicting the fact $u_{\infty}(x)=x_n+o(|x|)$. Up to a subsequence, $\frac{x_k}{|x_k|}\rightarrow e$ for some unit vector $e$. Up to a subsequence again, $$\frac{u(r_kx)}{r_k}\rightarrow \ \mbox{some blowdown function}\ V(x).$$ It is not difficult to see that $V(e)\geq 1$. Recalling that $V$ is weakly spacelike and $V(te_n)=t$ for $t\leq 0$, we have $$2=1-(-1)\leq V(e)-V(-e_n)\leq |e-(-e_n)|\leq 2,$$ implying $e=e_n$ and $V(e_n)=1$. By weakly spacelikeness again, we know $V(te_n)=t$ for $t\in [0,1]$. So $0$ can not be an isolated singularity of $V$ by Ecker's theorem. Thus
$V$ is an entire area maximizing hypersurface (and hence is a weakly spacelike hyperplane) and can only be the function $V(x)=x_n$.

Since we assumed $0\in \partial A$ and $A\subset \bar{B}_1$, one can verify that $u(x)\leq x_n+1$ by following the same way as in the proof of Theorem 1.1.

In the case that $M$ contains a future-directed light ray, we can apply the same argument.
\end{proof}

\begin{re}
In the case (b) of the above proof, one can also verify that $u(x)-x_n$ attains its maximum on $\partial A$. But the statement that $u(x)-x_n$ also attains its minimum on $\partial A$ can not be true in general. The example $w(x)$ in Section 3.2 provide a counterexample for this. Note that $w(te_n)<t-1$ for $t>1$ since $w$ is smooth and solves (1.1) in $\mathbb{R}^n\backslash \{te_n:t\in (-\infty,0]\}$. So if we let $A=\{0, e_n\}$ and consider $w$ as an exterior area maximizing function on
$\mathbb{R}^n\backslash \{0, e_n\}$, then $w(x)-x_n$ does not attain its minimum on $\partial A=\{0, e_n\}$. So far we don't know whether $w(x)-x_n$ is bounded from below.
\end{re}

\subsection{Exterior Dirichlet problems}

\begin{proof}[Proof of Theorem 1.3]
(i) and (ii) are symmetric. So we only prove (ii). Let $A$, $g$, $\psi$ (a weakly spacelike extension of $g$ into $\mathbb{R}^n\backslash A$) and $x_0$ be given as in the theorem. For simplicity, we assume $x_0=0$. Denote $\max\limits_{\partial A}(g(x)+|x|):=c^+$ and $\min\limits_{\partial A}(g(x)+|x|):=c^-$, then $$c^--|x|\leq g(x)\leq c^+-|x|\ \ \mbox{on}\ \partial A.$$ And it is easy to see that $\psi(x)\geq c^--|x|$ in $\mathbb{R}^n\backslash A$ by weakly spacelikeness of $\psi$. Define $\Psi(x):=\min(\psi(x), c^+-|x|)$. One can verify that $\Psi(x)$ is also a weakly spacelike extension of $g$ into $\mathbb{R}^n\backslash A$ and
$$c^--|x|\leq \Psi(x)\leq c^+-|x|\ \ \mbox{in}\ \mathbb{R}^n\backslash A.$$
Suppose $A\subset B_1$ without loss of generality. For $k=2,3,\cdots$, let $u_k$ be the area maximizing function in $B_k\backslash A$ satisfying $u_k=g$ on $\partial A$ and $u_k=\Psi$ on $\partial B_k$. The existence of such $u_k$ is due to Theorem 2.1. By Arzela-Ascoli theorem, up to a subsequence, $$u_k(x)\rightarrow u(x) \ \ \mbox{locally uniformly in}\ \mathbb{R}^n\backslash A.$$ By Lemma 2.1, $u(x)$ is area maximizing in $\mathbb{R}^n\backslash A$. By continuity, $u=g$ on $\partial A$. By comparison principle, $$c^--|x|\leq u(x)\leq c^+-|x|\ \ \mbox{in}\ \mathbb{R}^n\backslash A.$$

Now we prove (iii). Let $A$, $g$, $\psi$ and $a$ be given as in the theorem. For simplicity, we assume $a=e_n$. Denote $\max\limits_{\partial A}(g(x)-x_n):=c^+$ and $\min\limits_{\partial A}(g(x)-x_n):=c^-$, then $$c^-+x_n\leq g(x)\leq c^++x_n\ \ \mbox{on}\ \partial A.$$ Define $\Psi(x):=\max(\min(\psi(x), c^++x_n),c^-+x_n)$. One can still verify that $\Psi(x)$ is a weakly spacelike extension of $g$ into $\mathbb{R}^n\backslash A$ and
$$c^-+x_n\leq \Psi(x)\leq c^++x_n\ \ \mbox{in}\ \mathbb{R}^n\backslash A.$$ The rest things are totally same (replacing $-|x|$ by $x_n$) as in the proof of (ii) and we omit them.
\end{proof}

\begin{re}
As I mentioned in the introduction, unlike maximal hypersurfaces, we have no ideas about the uniqueness of solutions to the exterior Dirichlet problems for area maximizing hypersurfaces in all the three cases. However, thanks to the example in Section 3.2, we can still say something on this issue.

Let $A=\{0, e_n\}$, $g(x)=0$ on $\partial A=\{0, e_n\}$ and $a=e_n$, then $w(x)$ and $\tilde{w}(x)=\tilde{w}(x',x_n):=-w(x',1-x_n)$ are two different area maximizing functions in $\mathbb{R}^n\backslash A$
satisfying
\begin{equation}
u(x)=0\ \ \mbox{on}\ \partial A
\end{equation}
and
\begin{equation}
u(x)=x_n+o(|x|)\ \ \mbox{as}\ |x|\rightarrow \infty.
\end{equation}
By Remark 4.1, neither $w(x)$ nor $\tilde{w}(x)$ satisfies the condition
\begin{equation}
u(x)-x_n\ \mbox{attains its maximum and minimum on}\ \partial A.
\end{equation}
By Theorem 1.3, there exists an area maximizing function $\breve{w}$ in $\mathbb{R}^n\backslash A$ satisfying (4.1) and (4.3). That is to say, in this case, we have at least three different solutions to the exterior Dirichlet problem with the boundary conditions (4.1) and (4.2).

\end{re}

\section*{Acknowledgements}
The author would like to thank Professor Yu Yuan for many helpful discussions on this work. This paper was completed during the author's visit to University of Washington (Seattle). His visit was funded by China Scholarship Council. He would also like to thank Professor Yu Yuan for the invitation and to the Department of Mathematics for warm hospitality.

\bibliographystyle{elsarticle-num}

\end{document}